\documentclass[conference]{IEEEtran}
\usepackage{graphicx}  % Required for inserting images
\usepackage{amsmath}   % For math symbols and equations
\usepackage{cite}      % For citations

\begin{document}

\title{Analysis of Wind Power Integration in Electricity Markets  LMP Pricing}

\author{
	\IEEEauthorblockN{Narin Nemati{\dag}}
	\IEEEauthorblockA{
		Department of Applied Mathematics\\ Faculty of Mathematical Sciences\\ 
		University of Guilan, Rasht, Guilan, Iran\\
		Email: f.nematii78@gmail.com
	}
	\and
	\IEEEauthorblockN{Amin Ahmadi Kasani}
	\IEEEauthorblockA{
		Department of Mathematics, Statistics and Computer\\
		Science, College of Science, University of Tehran\\
		Email: aakasani@ut.ac.ir
	}
}

\maketitle
%\footnote{{\dag}corresponding author}
%\Author{Farshid Mehrdoust$^\dag$,  Seyedeh Fatemeh Nemati$^{\dag}$}

%\Address{$^{\dag}$Department of Applied Mathematics, Faculty of Mathematical Sciences\\ 
%University of Guilan, Rasht, Guilan, Iran}
        
%\Email{fmehrdoust@guilan.ac.ir, f.nematii78@gmail.com}
%\Markboth{F. Mehrdoust, S. F. Nemati}{Analysis of LMP pricing ...}

\begin{abstract}
	Wind energy has emerged as one of the most vital and economically viable forms of renewable energy. The integration of wind energy sources into power grids across the globe has been increasing substantially, largely due to the higher levels of uncertainty associated with wind energy compared to other renewable energy sources. This study focuses on analyzing the Locational Marginal Pricing (LMP) market model, with particular emphasis on the integration of wind power plants into substations. Furthermore, it examines a two-stage stochastic model for electricity markets employing LMP pricing, utilizing the Optimal Power Flow (OPF) method for the analysis.
\end{abstract}

%\AMS{2020}{91G15, 49M99.}
 
\section{Introduction}
In recent years, there has been a growing global focus on the utilization of renewable energy sources for electricity generation. The primary drivers of this shift are the limitations associated with fossil fuels and the environmental issues they exacerbate. Among the various forms of renewable energy, wind power plants have experienced the most significant growth worldwide. However, due to the inherent variability in wind power production, direct participation of these plants in the electricity market often fails to generate sufficient revenue, thereby discouraging investment in this sector. Electricity markets based on Locational Marginal Pricing (LMP), regarded as the most comprehensive market model, have been implemented in several countries across the globe. LMP is a mechanism used in electricity markets to determine the cost of electricity at different locations on the grid. In the context of wind power and electricity markets, LMP plays a critical role in addressing the fluctuations and unpredictability associated with wind power generation.\\

Recent research has focused on incorporating uncertainty factors into Locational Marginal Pricing (LMP) to strike an optimal balance between economic efficiency and system reliability when managing uncertainties in wind power output. A key element of this research involves integrating uncertainty components into LMP to price wind power and load uncertainties. These components, such as transmission line overload uncertainty costs and generation violation uncertainty costs, are incorporated into the LMP calculation to account for the financial impacts associated with uncertainties in wind power generation and demand at different locations. Additionally, the research underscores the necessity of balancing economic efficiency with system reliability when addressing the pricing of uncertainties related to wind power. Studies have also examined the application of LMP in multi-energy systems that rely extensively on renewable energy sources. These models aim to optimize economic dispatch by minimizing the additional costs associated with adjusting multi-energy units, thereby improving the grid's ability to adapt to variable wind power and complex pricing structures within electricity markets. In summary, the use of LMP in wind and electricity markets is primarily focused on addressing the challenges posed by the inherent uncertainty of wind power generation.\\
 
Incorporating uncertainty factors into Locational Marginal Pricing (LMP) calculations, alongside optimizing economic dispatch models, is a key research focus aimed at balancing economic efficiency with system reliability as renewable energy integration advances. LMP serves as one of the most widely utilized mechanisms for managing congestion in electricity markets \cite{litvinov2010design}. In \cite{liu2009locational}, several DC and AC optimal power flow (OPF) models are introduced to derive LMP values within wholesale electricity markets. Additionally, \cite{ahmed2011analysis} presents a stochastic model for LMP markets, specifically addressing the uncertainties associated with wind power generation. This model accounts for constraints related to power plant unit placement, grid limitations, and wind energy variability to assess its effects on market settlement, load distribution, and reserve requirements.\\

Furthermore, \cite{botterud2011wind} proposes a novel model for optimal wind energy trading in day-ahead (DA) electricity markets, under conditions of uncertainty regarding wind resources and prices. The work presented in \cite{morales2010simulating} investigates the impact of wind energy production on locational marginal prices in a competitive market setting. This analysis establishes a structural relationship between wind energy production and LMP values. The findings from this study provide a statistical framework for describing LMP fluctuations as a function of wind energy data and the market's structural characteristics.\\ 

Locational marginal prices (LMPs) are derived through an optimization method aimed at maximizing public welfare. Given the assumption of zero price elasticity of demand, cost minimization becomes equivalent to welfare maximization. As a result, the optimization approach employed in this study focuses on minimizing production costs. Additionally, it is assumed that perfect competition exists in the May network, meaning that the offers provided by electricity producers are set equal to their marginal costs. To ensure the economic and safe operation of power systems, the optimal power flow (OPF) method is typically utilized. In this study, the DC load flow method was applied, and LMPs were determined through the implementation of the OPF model. Mathematically, the LMP at each bus corresponds to the Lagrange multiplier associated with the power balance equation at that bus. As the penetration level of wind energy sources increases, the demand for both spinning and non-spinning reserves will also rise, thereby increasing operating costs. \\

In this study, the proposed method is utilized to compute the optimal reserve quantity and the associated costs of providing these reserves over a daily time horizon. The problem is formulated as a two-stage stochastic programming model. The first stage represents the electricity market, incorporating its constraints, while the second stage captures the operational and physical limitations of the power system. To evaluate the impact of wind energy's stochastic behavior, no other random processes, aside from wind energy, have been considered. For the sake of simplicity, it is assumed that the wind power plant is connected to a single grid bus. The structure of the paper is as follows: Section \ref{Sec2} presents the methodology for modeling wind energy uncertainty and the framework of the electricity market settlement model with locational marginal pricing. In Section \ref{Sec3}, the proposed model is applied to the 24-bus RTS network, and numerical results are presented.\\

\section{Formulation of the problem}\label{Sec2}
\subsection{Uncertainty modeling of wind resources}
In this section, wind energy production in the time period
$t$
is modeled as a random variable denoted by 
$P_{t}^{WP}$, $P^{WP}$ 
a set of dependent random 
$P_{t}^{WP}$
during the 
$N_{T}$
study time period, namely

\begin{center}
$P^{WP}=\{P_{t}^{WP}, t=1,2,...,N_{T}\}$
\end{center}

Constitutes a random process. In stochastic programing of the behavior of a stochastic process in the optimization problem, it is usually necessary to consider a set of scenarios 
$\omega$
that considers all possible states of the problem. In this article, each scenario 
$\omega \in \Omega$
is a wind energy production vector with
$N_{T}$
component, that is
$\omega=\{P_{1\omega}^{WP},...,P_{N_{T}\omega}^{WP}\}$
scenario. In other words, 
$P_{t\omega}^{WP}$
represents the realization of the random variable 
$P_{t}^{WP}$
in the scenario 
$\omega$. Each scenario $\omega$  
has a probability of occurrence 
$\pi_{\omega}$,
so that the sum of the probabilities of all scenarios is equal to, that is

\begin{equation}
\displaystyle\sum\limits_{\omega=1}^{N_{\Omega}} \pi_{\omega}= 1.
\end{equation}

\subsection{LMP market stochastic model}
The objective function of the problem that must be minimized is the cost function as follows 

\begin{center}
$EC=\displaystyle\sum\limits_{t=1}^{N_{T}}EC_{t}=\displaystyle\sum\limits_{t=1}^{N_{T}}\displaystyle\sum\limits_{i=1}^{N_{G}}C_{it}^{SU}$
\end{center}

\begin{center}
$+\displaystyle\sum\limits_{t=1}^{N_{T}}d_{t}[\displaystyle\sum\limits_{i=1}^{N_{G}}\displaystyle\sum\limits_{m=1}^{N_{it}}\lambda_{Git}(m)p_{Git}(m)-\displaystyle\sum\limits_{j=1}^{N_{L}}\lambda_{Ljt}L_{jt}^{S}$
\end{center}

\begin{center}
$+\displaystyle\sum\limits_{i=1}^{N_{G}}(C_{it}^{R^{U}}R_{it}^{U}+C_{it}^{R^{D}}R_{it}^{D}+C_{it}^{R^{NS}}R_{it}^{NS})
+\displaystyle\sum\limits_{j=1}^{N_{L}}(C_{jt}^{R^{U}}R_{jt}^{U}+C_{jt}^{R^{D}}R_{jt}^{D})+\lambda_{t}^{WP}P_{t}^{WP,S}]$
\end{center}

\begin{center}
$+\displaystyle\sum\limits_{\omega=1}^{N_{\omega}}\pi_{\omega}\{{\displaystyle\sum\limits_{t=1}^{N_{T}}\displaystyle\sum\limits_{i=1}^{N_{G}}C_{it\omega}}+\displaystyle\sum\limits_{t=1}^{N_{T}}d_{t}[\displaystyle\sum\limits_{i=1}^{N_{G}}\displaystyle\sum\limits_{m=1}^{N_{it}}\lambda_{Git}(m)r_{Git\omega}(m)$
\end{center}
\begin{equation}
+\displaystyle\sum\limits_{j=1}^{N_{L}}\lambda_{Ljt}(r_{jt\omega}^{U}-r_{jt\omega}^{D})+\displaystyle\sum\limits_{j=1}^{N_{L}}V_{jt}L_{jt\omega}]\}.
\end{equation}

In this regard, $EC_{t}$ represents the total system cost during the time period $t$, while $C_{it}^{R^{U}}$, $C_{it}^{R^{D}}$, and $C_{it}^{R^{NS}}$ denote the proposed costs for incremental, decremental, and non-spinning reserves for unit $i$ during time period $t$, respectively. Similarly, $C_{jt}^{R^{U}}$ and $C_{jt}^{R^{D}}$ represent the proposed costs for incremental and spinning reserves for load $j$ during time period $t$, respectively.The objective function consists of nine components, which are explained as follows: 1-Proposed cost of starting up power plant units, 2-Proposed cost of energy production from power plant units minus consumer utility, 3-Proposed cost of providing incremental/decremental spinning reserves and non-spinning reserves from power plant units, 4-Proposed cost of providing incremental/decremental spinning reserves and non-spinning reserves from the load, 5-Proposed cost of wind power generation, 6-Costs associated with the start-up and shutdown scheduling of power plant units, 7-Costs related to the effective use of incremental/decremental spinning reserves and non-spinning reserves by power plant units, 8-Costs related to the effective use of incremental/decremental spinning reserves by the load, 9-Load shedding costs.In this context, it is assumed that wind power plants are not competitive entities. Therefore, the following problem is considered in the objective function.\\
\begin{center}
 $\lambda_{t}^{WP}=0.$
\end{center}

\subsection{Limitations of the Problem}

\subsubsection{Constraints Related to the Electricity Market and First-Stage Variables (Independent of the $\omega$ Scenario)}
The following limitations are associated with the electricity market, specifically focusing on first-stage variables that are not dependent on the $\omega$ scenario:

\begin{enumerate}
	\item Market Balance Constraint
	
	The electricity market must maintain balance between supply and demand at each time step. This is represented by the following equation:
	
	\begin{equation}
		\displaystyle\sum\limits_{i=1}^{N_{G}}P_{it}^{s} + P_{t}^{WP,s} = \displaystyle\sum\limits_{j=1}^{N_{L}}L_{jt}^{s}, \quad \forall t
	\end{equation}
	
	Where $P_{it}^{s}$ denotes the power generated by generator $i$ at time $t$ in scenario $s$, $P_{t}^{WP,s}$ represents the wind power generation at time $t$ in scenario $s$, and $L_{jt}^{s}$ refers to the load demand at bus $j$ at time $t$ in scenario $s$. This equation ensures that the total generation (including wind power) matches the total load demand across the system for each time period.

\item Generation Constraints

The following constraints apply to the power generation limits of the generating units:

\begin{equation}
	P_{i}^{min}u_{it} \leq P_{it}^{s} \leq P_{i}^{max} u_{it}, \quad \forall i,t
\end{equation}

This equation ensures that the power output of unit $i$ at time $t$ remains within its minimum ($P_{i}^{min}$) and maximum ($P_{i}^{max}$) generation limits, depending on whether the unit is operational ($u_{it}$).

\begin{equation}
	0 \leq p_{Git}(m) \leq p_{Git}^{max}(m), \quad \forall m,i,t
\end{equation}

This constraint specifies that the output of each generator $m$ within unit $i$ at time $t$ must be non-negative and bounded by its maximum output ($p_{Git}^{max}(m)$).

\begin{equation}
	P_{it}^{s} = \displaystyle\sum\limits_{m=1}^{N_{it}} p_{Git}(m),\quad \forall i,t
\end{equation}

Here, the total power output of unit $i$ at time $t$ ($P_{it}^{s}$) is the sum of the outputs of all $N_{it}$ individual generators.

In this context, $P_{i}^{max}$ and $P_{i}^{min}$ represent the maximum and minimum power output capacities of generating unit $i$, respectively.

\item Wind Power Plant Production Constraints

The following constraints apply to wind power generation:

\begin{equation}
	P_{t}^{WP,min} \leq P_{t}^{WP,s} \leq P_{t}^{WP,max}, \quad \forall t
\end{equation}

This equation ensures that the wind power generated at time $t$ ($P_{t}^{WP,s}$) remains within the minimum ($P_{t}^{WP,min}$) and maximum ($P_{t}^{WP,max}$) production limits.

In this context, $P_{t}^{WP,min}$ and $P_{t}^{WP,max}$ represent the lower and upper bounds of the wind power generation capacity, as specified in the wind power proposal.

\item Load Constraints

The load demand at each bus is subject to the following constraints:

\begin{equation}
	L_{jt}^{s,min} \leq L_{jt}^{s} \leq L_{jt}^{s,max}, \quad \forall j,t
\end{equation}

Here, $L_{jt}^{s,min}$ and $L_{jt}^{s,max}$ represent the minimum and maximum load demand at bus $j$ during time period $t$, which are specified as part of the demand-side offer. In cases where the load is price-inelastic, the upper and lower limits are equal, implying that $L_{jt}^{s,min} = L_{jt}^{s} = L_{jt}^{s,max}$.

\item Scheduled Reserve Determination Constraints

\textbf{5-1. Supply Side:}

The scheduled reserve constraints on the supply side are defined as follows:

\begin{equation}
	0 \leq R_{it}^{U} \leq R_{it}^{U,max}u_{it}, \quad \forall i,t
\end{equation}

\begin{equation}
	0 \leq R_{it}^{D} \leq R_{it}^{D,max}u_{it}, \quad \forall i,t
\end{equation}

\begin{equation}
	0 \leq R_{it}^{NS} \leq R_{it}^{NS,max}(1-u_{it}), \quad \forall i,t
\end{equation}

In these equations, $R_{it}^{U}$ and $R_{it}^{D}$ represent the ramp-up and ramp-down reserves for generating unit $i$ during time period $t$, while $R_{it}^{NS}$ denotes the non-spinning reserve. The variables $R_{it}^{U,max}$, $R_{it}^{D,max}$, and $R_{it}^{NS,max}$ define the maximum ramp-up, ramp-down, and non-spinning reserve limits, respectively. The binary variable $u_{it}$ indicates whether the generating unit is operational.

\textbf{5-2. Demand Side:}

The reserve constraints on the demand side are defined as follows:

\begin{equation}
	0 \leq R_{jt}^{U} \leq R_{jt}^{U,max}, \quad \forall j,t
\end{equation}

\begin{equation}
	0 \leq R_{jt}^{D} \leq R_{jt}^{D,max}, \quad \forall j,t
\end{equation}

Here, $R_{jt}^{U}$ and $R_{jt}^{D}$ denote the upward and downward reserves for load $j$ during time period $t$. The variables $R_{jt}^{U,max}$ and $R_{jt}^{D,max}$ represent the corresponding maximum reserve limits.

\item Start-Up Cost Constraints

The following constraints govern the start-up costs for generating units:

\begin{equation}
	C_{it}^{SU} \geq \lambda_{it}^{SU}(u_{it}-u_{i,t-1}), \quad \forall i,t
\end{equation}

\begin{equation}
	C_{it}^{SU} \geq 0, \quad \forall i,t
\end{equation}

In these equations, $C_{it}^{SU}$ represents the start-up cost for generating unit $i$ at time $t$, while $\lambda_{it}^{SU}$ denotes the start-up cost coefficient. The term $(u_{it} - u_{i,t-1})$ captures the binary change in the unit’s status, where $u_{it}$ is the unit's status at time $t$ and $u_{i,t-1}$ is the status at the previous time period.

\item Constraints Related to System Performance and Second-Stage Variables (Dependent on the $\omega$ Scenario)

These constraints focus on the system’s actual performance and take into account the variables in the second stage, which depend on the specific $\omega$ scenario. Further detail on these constraints is provided in subsequent sections.\\
\end{enumerate}
\begin{enumerate}
	\item Power Balance Constraints
	
	\textbf{1-1. Power balance at each node $n$ (excluding the node $n^{\prime}$ connected to the wind power plant):}
	
	The power balance at node $n$ for nodes not connected to the wind power plant is described by the following equation:
	
	\begin{multline}
		\displaystyle\sum\limits_{i:(i,n) \in M_{G}} P_{it\omega}^{G} - \displaystyle\sum\limits_{j:(j,n) \in M_{L}} (L_{jt\omega} ^{c} -  L_{jt\omega}) - \displaystyle\sum\limits_{r:(n,r)\in \Lambda} \\ f_{t\omega}(n,r) = 0, \quad \forall n \neq n^{\prime}, \quad \forall t, \omega
	\end{multline}
	
	In this equation, $P_{it\omega}^{G}$ represents the power generated by generator $i$ at time $t$ under scenario $\omega$, and $L_{jt\omega}^{c}$ and $L_{jt\omega}$ represent the curtailed and actual load demand, respectively. The term $f_{t\omega}(n,r)$ is the power flow between nodes $n$ and $r$, and $\Lambda$ denotes the set of transmission lines.
	
	\textbf{1-2. Power balance at node $n^{\prime}$ (connected to the wind power plant):}
	
	For the node $n^{\prime}$ connected to the wind power plant, the power balance is adjusted to account for wind power injection:
	
	\begin{multline}
		\displaystyle\sum\limits_{i:(i,n) \in M_{G}} P_{it\omega}^{G} - \displaystyle\sum\limits_{j:(j,n) \in M_{L}} (L_{jt\omega} ^{c} \\ -  L_{jt\omega}) + P_{t\omega}^{WP} - \displaystyle\sum\limits_{r:(n,r)\in \Lambda} f_{t\omega}(n,r) = 0, \quad \\
		 \forall n = n^{\prime}, \quad \forall t, \omega
	\end{multline}
	
	Here, $P_{t\omega}^{WP}$ represents the wind power injected into node $n^{\prime}$ at time $t$ under scenario $\omega$.
	
	\textbf{1-3. Power flow through transmission line $(n,r)$:}
	
	The flow of power along the transmission line connecting nodes $n$ and $r$ is given by:
	
	\begin{multline}
		f_{t\omega}(n,r) = \frac{P_{t\omega}^{loss}(n,r)}{2} + \\ B(n,r)(\delta_{nt\omega} - \delta_{rt\omega}), \quad \forall(n,r) \in \Lambda, \quad \forall t, \omega
	\end{multline}
	
	In this equation, $P_{t\omega}^{loss}(n,r)$ denotes the power losses along the line $(n,r)$, and $B(n,r)$ represents the susceptance of the transmission line. The terms $\delta_{nt\omega}$ and $\delta_{rt\omega}$ are the voltage phase angles at nodes $n$ and $r$, respectively.
	
	\item Generation Constraints
	
	The generation capacity of each unit is constrained as follows:
	
	\begin{equation}
		P_{it\omega}^{G} \geq P_{i}^{min} v_{it\omega}, \quad \forall i, t, \omega
	\end{equation}
	
	\begin{equation}
		P_{it\omega}^{G} \leq P_{i}^{max} v_{it\omega}, \quad \forall i, t, \omega
	\end{equation}
	
	Here, $P_{i}^{min}$ and $P_{i}^{max}$ represent the minimum and maximum generation limits for generator $i$, and $v_{it\omega}$ is a binary variable indicating whether the unit is operational at time $t$ under scenario $\omega$.
	
	\item Transmission Line Capacity Constraints
	
	The transmission capacity of each line is limited as follows:
	
	\begin{equation}
		-f^{max}(n,r) \leq f_{t\omega}(n,r) \leq f^{max}(n,r), \quad \forall(n,r) \in \Lambda
	\end{equation}
	
	This equation ensures that the power flow through the transmission line $(n,r)$ remains within the maximum allowable limits, $f^{max}(n,r)$.
\end{enumerate}

\section{Numerical experiments}\label{Sec3}
\subsection{Sample network}
\begin{figure}[h] 
 \centering{\includegraphics[width=9cm,height=10cm ]{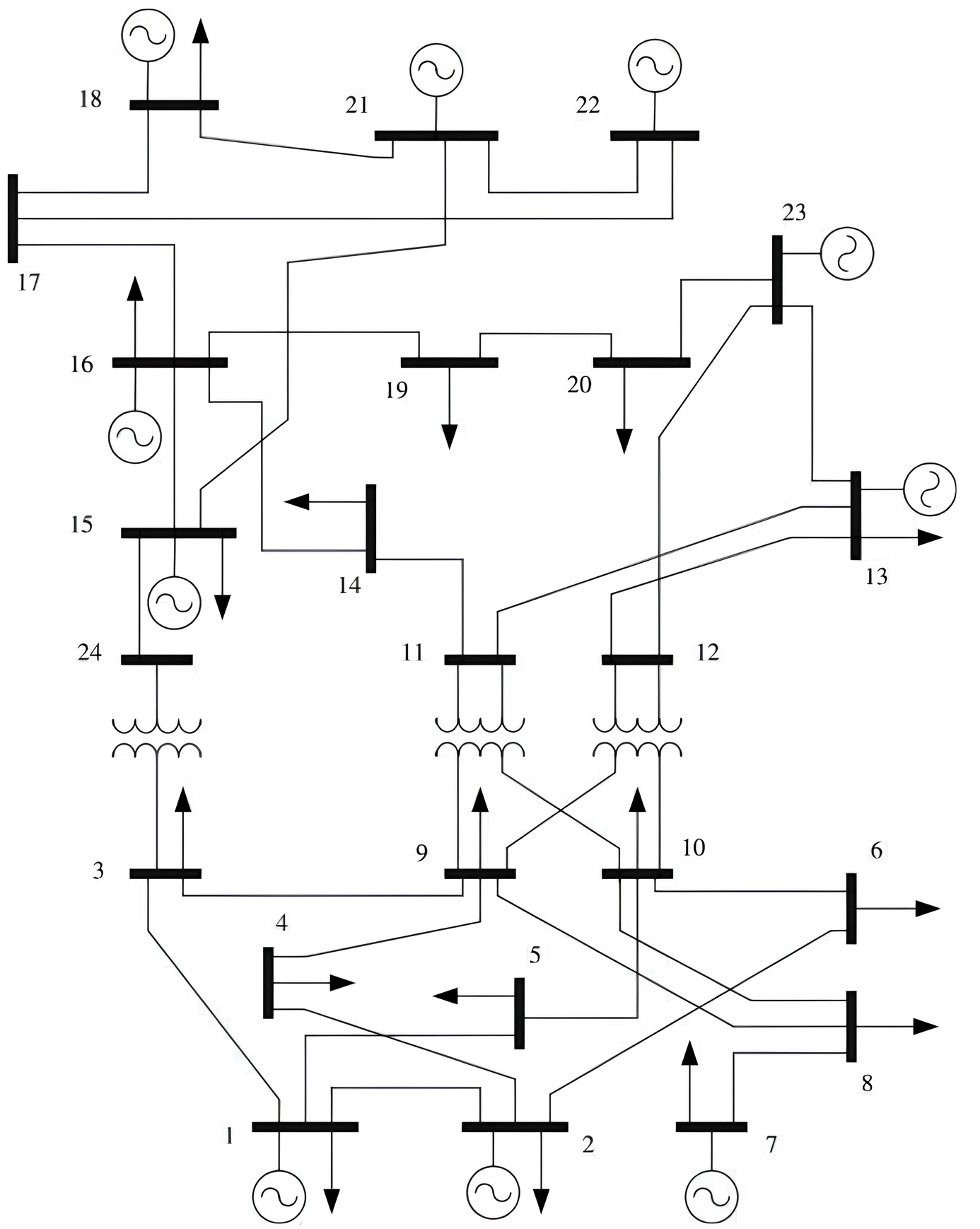}} 
\caption[ ]
{ RTS 24 bus network \cite{farashbashi2015comparative,ordoudis2016updated}}
 \label{fig:1}
\end{figure}

The network used for evaluating the market settlement model with locational marginal pricing (LMP) in this study is based on the 24-bus IEEE RTS (Reliability Test System) network and its single-zone configuration. As illustrated in Figure \ref{fig:1}, the network consists of 34 transmission lines. The resistance of these lines is assumed to be zero, meaning real power losses are neglected in the analysis.

 Each power plant in the system submits offers for the maximum possible values of its upward and downward spinning reserves, defined as the difference between its maximum and minimum generation capacities $(p^{max}-p^{min})$. These reserves are offered at a cost equivalent to 25\% of the highest energy production cost bid. However, only the power plants connected to buses 7, 15, and 16 offer non-spinning reserves. These plants also provide the maximum possible amount of their non-spinning reserves $(p^{max})$, offered at a cost equivalent to 20\% of their highest energy production cost.
 
 All generators offer the same rate of 50 $\frac{\$}{MWh}$ for incremental and decremental spinning reserves. Furthermore, each load in the system has the flexibility to increase or decrease its consumption by up to 20\% of its planned demand at any hour to provide upward or downward spinning reserves. If the load exceeds these limits, the lost load is valued at 2000 $\frac{\$}{MWh}$.\\

\subsection{Wind energy scenario information }
In the analysis conducted for this project, wind energy scenarios, as outlined in Table \ref{table1}, were used. High and low forecast scenarios were considered over a 4-hour time period.\\

\begin{table}[h] 
\begin{center}
\caption{Wind energy scenarios} 

\label{table1}
\end{center}
%\begin{center}
\[
\begin{tabular}{|c|c|c|c|}
\hline 
$\text{Period}$ & $ \text{As forecast} $ & $\text{High}$ & $\text{Low}$ \\
\hline 
$t=1$ & $8$ & $11$ & $1$  \\ 
\hline 
$t=2 $ & $21$ & $28$ & $11$   \\ 
\hline 
$t=3$ & $36$ & $48$ & $27 $ \\ 
\hline 
$t=4$ & $9$ & $13$ & $7 $ \\ 
\hline
\end{tabular}
\]
%\end{center}
\end{table}
\subsection{Analysis of the results}
\begin{itemize}
\item[1-] Investigating the impact of wind resource uncertainty

The uncertainty of wind resources is defined as the difference between the expected scenario of energy production from wind and the upper and lower scenarios. The degree of uncertainty is expressed as a percentage from $0\%$ to $100\%$, where an uncertainty of $x\%$ means that the upper and lower scenarios are equal to the predicted wind energy production scenario multiplied by $(1+x\%)$ and $(1-x\%)$, respectively. In this section, the uncertainty of wind resources is increased by changing the probability of occurrence of the scenarios, and the results of the software implementation for the three modes are presented in Table \ref{table2}.\\
\begin{table}[h] 
\begin{center}
\caption{Effect of uncertainty on the objective function of the problem} 

\label{table2}
\end{center}
\begin{center}
\resizebox{10cm}{!}{%
{\begin{tabular}{|c|c|c|c|}
\hline 
$ $ & $\text{Percentage of uncertainty}(40\%) $ & $\text{Percentage of uncertainty}(50\%) $ & $\text{Percentage of uncertainty}(60\%) $ \\
\hline 
$\text{Objective function}$ & $40528$ & $41239$ & $42123$  \\ 
\hline 
\end{tabular}}%
}
\end{center}
\end{table}

As uncertainty increases, the optimal value of the market objective function rises. Since there is no cost associated with wind energy production in this project, an increase in uncertainty leads to a decrease in the probability of wind energy production. Consequently, other power plants increase their production, resulting in higher total costs. The results of the implementation of the proposed market settlement model, with uncertainty increasing from $40\%$ to $50\%$, have also been examined. According to the results, as uncertainty increases, power plant production decreases, which is compensated for by the reserve market. Considering that, according to the problem's assumptions, each power plant offers the maximum possible values for its incremental and decremental spinning reserves at a cost equal to $25\%$ of the highest energy production cost, an increase in LMP is expected at the bus to which this power plant is connected. \\

\item[2-] Investigating the effect of wind sources on line congestion conditions.
In this section, to create congestion, the capacity of the line connecting bus 1 to bus 2 has been reduced, and its effect on LMP values in certain buses has been investigated. 
The optimal value of the market objective function has changed, as shown in Table \ref{table3}.\\

\begin{table}[h] 
\begin{center}
\caption{Effect of creating density on the objective function of the problem} 

\label{table3}
\end{center}
\begin{center}
{\begin{tabular}{|c|c|c|}
\hline 
$ $ & $f_{max}(1,2)=175MW$ & $f_{max}(1,2)=10MW$  \\
\hline 
$\text{Objective function}$ & $40528$ & $41896$   \\ 
\hline 
\end{tabular}}
\end{center}
\end{table} 

In this case, power station 2, which is connected to bus 2 and whose optimal production would be zero under no congestion, will experience an increase in production, leading to a rise in the LMP value at this bus compared to other nearby buses.\\
\begin{figure}[h] 
\centering{\includegraphics[width=10cm,height=8cm ]{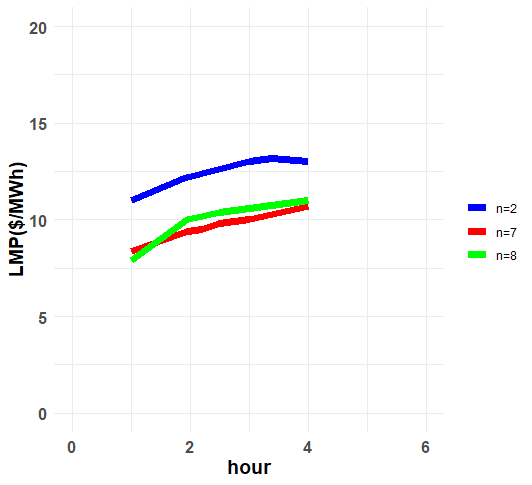}} 
\caption[ ]
{ Comparison of LMP values in adjacent tires with tire 2 in compression mode }
 \label{fig:4}
\end{figure} 

In this section, the effect of increasing the penetration coefficient of the wind source connected to bus 2 on the LMP values at this bus is investigated. According to the figure, with a $5\%$ increase in the penetration coefficient of the wind source, the LMP values decreased. This decrease is due to the reduction in congestion on the connecting line between buses 1 and 2 caused by the wind source.\\
  
\begin{figure}[h] 
\centering{\includegraphics[width=10cm,height=8cm ]{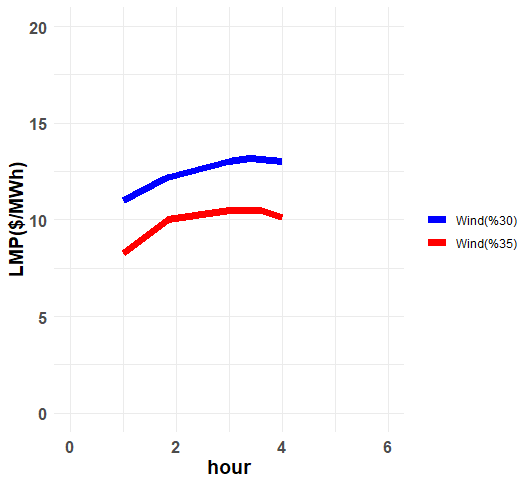}} 
\caption[ ]
{ Effect of increasing the penetration coefficient of the wind source connected to tire 2 on the LMP values of this tire  }
 \label{fig:5}
\end{figure}
\item[3-] Investigating the effect of increasing the penetration coefficient of wind sources on   LMP values in the non-condensation mode

In this section, the capacity of the wind resource connected to bus 2 has been increased, and its effect on the value of locational marginal prices at this bus has been investigated. The results obtained from the calculations are shown in Table \ref{table5}. In this analysis, the penetration coefficient of the wind resource is calculated as follows. \\
    \begin{table}[h] 
\begin{center}
\caption{considered penetration coefficient} 

\label{table5}
\end{center}
\begin{center}
\begin{tabular}{|c|c|}
\hline 
$ $ & $\text{Penetration coefficient}$  \\
\hline 
$W_{1}$ & $35\%$  \\ 
\hline 
$W_{2}$ & $40\%$ \\ 
\hline
$W_{3}$ & $45\%$ \\ 
\hline
\end{tabular}
\end{center}
\end{table}

\begin{figure}[h] 
\centering{\includegraphics[width=10cm,height=8cm ]{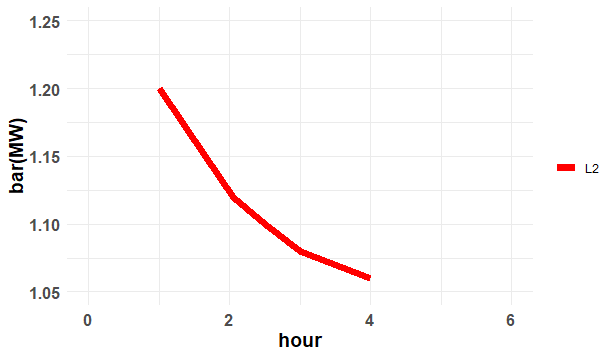}} 
\caption[ ]
{ load profile connected to bus 2 }
 \label{fig:7}
\end{figure}

\begin{figure}[h] 
\centering{\includegraphics[width=10cm,height=8cm ]{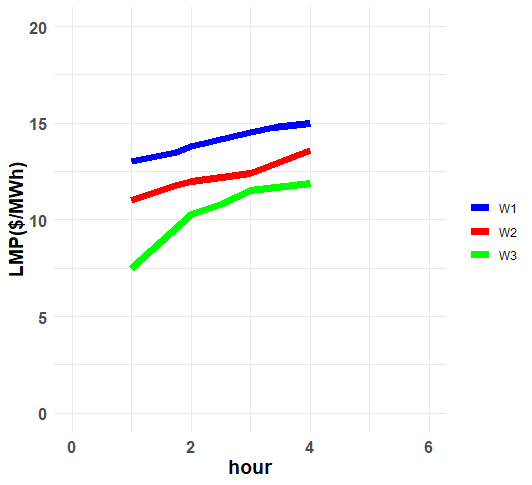}} 
\caption[ ]
{ Influence of the penetration coefficient of wind sources on LMP values of tire 2 in the non-condensation mode }
 \label{fig:8}
\end{figure}

According to the obtained graph, the LMP values at the bus connected to the wind power plant begin to decrease as the capacity of the wind power plant increases. This is due to the contribution of wind energy, which is provided at zero cost, to the injected energy at this bus. As a result, congestion and costs at this bus are reduced, leading to a lower locational marginal price. When the load profile is considered over a 4-hour period, the decrease in locational marginal prices is more significant during peak load hours. \\
\end{itemize}

\section{Conclusion}
Our analysis has introduced a new LMP pricing framework that takes into account the inherent randomness of wind energy in different regions of the power system. Using sophisticated statistical models, we can predict wind power generation at different grid locations. This allows for real-time adjustment of LMPs and reflects the impact of wind on supply and transmission constraints. This innovation fosters a more efficient market, optimizes power dispatch decisions, and promotes investment in energy storage to address the intermittent nature of wind.\\

In this paper, the electricity market with LMP pricing is described as a two-stage stochastic programming problem. The first stage represents the electricity market and its related constraints, while the second stage represents the power system, its operation, and physical limitations. With the analysis performed on the sample network, the following results have been obtained:\\
 
\begin{enumerate}
	\item The use of wind energy sources under conditions of high uncertainty increases the expected cost of the network. In other words, uncertainty reduces the impact of wind resources on lowering the overall system cost to some extent.
	\item Congestion in the network lines causes the LMP at buses connected to congested lines to increase. However, this can be mitigated by connecting wind power plants to these buses, which will result in a decrease in LMP.
	\item In the absence of congestion in the network lines, an increase in the penetration coefficient of wind resources will lead to a decrease in LMP values. This decrease will be more significant during peak hours.
\end{enumerate}

\section{Author Contributions}
N.N. conceptualized the study, conducted data analysis and interpretation, and drafted the manuscript.\\
A.A.K. contributed to reviewing and editing the final draft and was responsible for overall supervision.\\

\bibliographystyle{IEEEtran}
\bibliography{citation}

\begin{thebibliography}{1}
\providecommand{\url}[1]{#1}
\csname url@samestyle\endcsname
\providecommand{\newblock}{\relax}
\providecommand{\bibinfo}[2]{#2}
\providecommand{\BIBentrySTDinterwordspacing}{\spaceskip=0pt\relax}
\providecommand{\BIBentryALTinterwordstretchfactor}{4}
\providecommand{\BIBentryALTinterwordspacing}{\spaceskip=\fontdimen2\font plus
\BIBentryALTinterwordstretchfactor\fontdimen3\font minus
  \fontdimen4\font\relax}
\providecommand{\BIBforeignlanguage}[2]{{%
\expandafter\ifx\csname l@#1\endcsname\relax
\typeout{** WARNING: IEEEtran.bst: No hyphenation pattern has been}%
\typeout{** loaded for the language `#1'. Using the pattern for}%
\typeout{** the default language instead.}%
\else
\language=\csname l@#1\endcsname
\fi
#2}}
\providecommand{\BIBdecl}{\relax}
\BIBdecl

\bibitem{litvinov2010design}
E.~Litvinov, ``Design and operation of the locational marginal prices-based
  electricity markets,'' \emph{IET generation, transmission \& distribution},
  vol.~4, no.~2, pp. 315--323, 2010.

\bibitem{liu2009locational}
H.~Liu, L.~Tesfatsion, and A.~Chowdhury, ``Locational marginal pricing basics
  for restructured wholesale power markets,'' in \emph{2009 IEEE Power \&
  Energy Society General Meeting}.\hskip 1em plus 0.5em minus 0.4em\relax IEEE,
  2009, pp. 1--8.

\bibitem{ahmed2011analysis}
M.~H. Ahmed, K.~Bhattacharya, and M.~Salama, ``Analysis of uncertainty model to
  incorporate wind penetration in lmp-based energy markets,'' in \emph{2011 2nd
  International Conference on Electric Power and Energy Conversion Systems
  (EPECS)}.\hskip 1em plus 0.5em minus 0.4em\relax IEEE, 2011, pp. 1--8.

\bibitem{botterud2011wind}
A.~Botterud, Z.~Zhou, J.~Wang, R.~J. Bessa, H.~Keko, J.~Sumaili, and
  V.~Miranda, ``Wind power trading under uncertainty in lmp markets,''
  \emph{IEEE Transactions on power systems}, vol.~27, no.~2, pp. 894--903,
  2011.

\bibitem{morales2010simulating}
J.~M. Morales, A.~J. Conejo, and J.~Perez-Ruiz, ``Simulating the impact of wind
  production on locational marginal prices,'' \emph{IEEE Transactions on Power
  Systems}, vol.~26, no.~2, pp. 820--828, 2010.

\bibitem{farashbashi2015comparative}
S.-M. Farashbashi-Astaneh, W.~Hu, and Z.~Chen, ``Comparative study between two
  market clearing schemes in wind dominant electricity markets,'' \emph{IET
  Generation, Transmission \& Distribution}, vol.~9, no.~15, pp. 2215--2223,
  2015.

\bibitem{ordoudis2016updated}
C.~Ordoudis, P.~Pinson, J.~M.~M. Gonz{\'a}lez, and M.~Zugno, ``An updated
  version of the ieee rts 24-bus system for electricity market and power system
  operation studies.'' 2016.

\end{thebibliography}

\end{document}